\input amstex
\input amsppt.sty   
%
\catcode`\@=11

\def\East#1#2{\setboxz@h{$\m@th\ssize\;{#1}\;\;$}%
 \setbox@ne\hbox{$\m@th\ssize\;{#2}\;\;$}\setbox\tw@\hbox{$\m@th#2$}%
 \dimen@\minaw@
 \ifdim\wdz@>\dimen@ \dimen@\wdz@ \fi  \ifdim\wd@ne>\dimen@ \dimen@\wd@ne \fi
 \ifdim\wd\tw@>\z@
  \mathrel{\mathop{\hbox to\dimen@{\rightarrowfill}}\limits^{#1}_{#2}}%
 \else
  \mathrel{\mathop{\hbox to\dimen@{\rightarrowfill}}\limits^{#1}}%
 \fi}
\def\West#1#2{\setboxz@h{$\m@th\ssize\;\;{#1}\;$}%
 \setbox@ne\hbox{$\m@th\ssize\;\;{#2}\;$}\setbox\tw@\hbox{$\m@th#2$}%
 \dimen@\minaw@
 \ifdim\wdz@>\dimen@ \dimen@\wdz@ \fi \ifdim\wd@ne>\dimen@ \dimen@\wd@ne \fi
 \ifdim\wd\tw@>\z@
  \mathrel{\mathop{\hbox to\dimen@{\leftarrowfill}}\limits^{#1}_{#2}}%
 \else
  \mathrel{\mathop{\hbox to\dimen@{\leftarrowfill}}\limits^{#1}}%
 \fi}
\font\arrow@i=lams1
\font\arrow@ii=lams2
\font\arrow@iii=lams3
\font\arrow@iv=lams4
\font\arrow@v=lams5
\newbox\zer@
\newdimen\standardcgap
\standardcgap=40\p@
\newdimen\hunit
\hunit=\tw@\p@
\newdimen\standardrgap
\standardrgap=32\p@
\newdimen\vunit
\vunit=1.6\p@
\def\Cgaps#1{\RIfM@
  \standardcgap=#1\standardcgap\relax \hunit=#1\hunit\relax
 \else \nonmatherr@\Cgaps \fi}
\def\Rgaps#1{\RIfM@
  \standardrgap=#1\standardrgap\relax \vunit=#1\vunit\relax
 \else \nonmatherr@\Rgaps \fi}
\newdimen\getdim@
\def\getcgap@#1{\ifcase#1\or\getdim@\z@\else\getdim@\standardcgap\fi}
\def\getrgap@#1{\ifcase#1\getdim@\z@\else\getdim@\standardrgap\fi}
\def\cgaps#1{\RIfM@
 \cgaps@{#1}\edef\getcgap@##1{\i@=##1\relax\the\toks@}\toks@{}\else
 \nonmatherr@\cgaps\fi}
\def\rgaps#1{\RIfM@
 \rgaps@{#1}\edef\getrgap@##1{\i@=##1\relax\the\toks@}\toks@{}\else
 \nonmatherr@\rgaps\fi}
\def\Gaps@@{\gaps@@}
\def\cgaps@#1{\toks@{\ifcase\i@\or\getdim@=\z@}%
 \gaps@@\standardcgap#1;\gaps@@\gaps@@
 \edef\next@{\the\toks@\noexpand\else\noexpand\getdim@\noexpand\standardcgap
  \noexpand\fi}%
 \toks@=\expandafter{\next@}}
\def\rgaps@#1{\toks@{\ifcase\i@\getdim@=\z@}%
 \gaps@@\standardrgap#1;\gaps@@\gaps@@
 \edef\next@{\the\toks@\noexpand\else\noexpand\getdim@\noexpand\standardrgap
  \noexpand\fi}%
 \toks@=\expandafter{\next@}}
\def\gaps@@#1#2;#3{\mgaps@#1#2\mgaps@
 \edef\next@{\the\toks@\noexpand\or\noexpand\getdim@
  \noexpand#1\the\mgapstoks@@}%
 \global\toks@=\expandafter{\next@}%
 \DN@{#3}%
 \ifx\next@\Gaps@@\gdef\next@##1\gaps@@{}\else
  \gdef\next@{\gaps@@#1#3}\fi\next@}
\def\mgaps@#1{\let\mgapsnext@#1\FN@\mgaps@@}
\def\mgaps@@{\ifx\next\space@\DN@. {\FN@\mgaps@@}\else
 \DN@.{\FN@\mgaps@@@}\fi\next@.}
\def\mgaps@@@{\ifx\next\w\let\next@\mgaps@@@@\else
 \let\next@\mgaps@@@@@\fi\next@}
\newtoks\mgapstoks@@
\def\mgaps@@@@@#1\mgaps@{\getdim@\mgapsnext@\getdim@#1\getdim@
 \edef\next@{\noexpand\getdim@\the\getdim@}%
 \mgapstoks@@=\expandafter{\next@}}
\def\mgaps@@@@\w#1#2\mgaps@{\mgaps@@@@@#2\mgaps@
 \setbox\zer@\hbox{$\m@th\hskip15\p@\tsize@#1$}%
 \dimen@\wd\zer@
 \ifdim\dimen@>\getdim@ \getdim@\dimen@ \fi
 \edef\next@{\noexpand\getdim@\the\getdim@}%
 \mgapstoks@@=\expandafter{\next@}}
\def\changewidth#1#2{\setbox\zer@\hbox{$\m@th#2}%
 \hbox to\wd\zer@{\hss$\m@th#1$\hss}}
\atdef@({\FN@\ARROW@}
\def\ARROW@{\ifx\next)\let\next@\OPTIONS@\else
 \DN@{\csname\string @(\endcsname}\fi\next@}
\newif\ifoptions@
\def\OPTIONS@){\ifoptions@\let\next@\relax\else
 \DN@{\options@true\begingroup\optioncodes@}\fi\next@}
\newif\ifN@
\newif\ifE@
\newif\ifNESW@
\newif\ifH@
\newif\ifV@
\newif\ifHshort@
\expandafter\def\csname\string @(\endcsname #1,#2){%
 \ifoptions@\let\next@\endgroup\else\let\next@\relax\fi\next@
 \N@false\E@false\H@false\V@false\Hshort@false
 \ifnum#1>\z@\E@true\fi
 \ifnum#1=\z@\V@true\tX@false\tY@false\a@false\fi
 \ifnum#2>\z@\N@true\fi
 \ifnum#2=\z@\H@true\tX@false\tY@false\a@false\ifshort@\Hshort@true\fi\fi
 \NESW@false
 \ifN@\ifE@\NESW@true\fi\else\ifE@\else\NESW@true\fi\fi
 \arrow@{#1}{#2}%
 \global\options@false
 \global\scount@\z@\global\tcount@\z@\global\arrcount@\z@
 \global\s@false\global\sxdimen@\z@\global\sydimen@\z@
 \global\tX@false\global\tXdimen@i\z@\global\tXdimen@ii\z@
 \global\tY@false\global\tYdimen@i\z@\global\tYdimen@ii\z@
 \global\a@false\global\exacount@\z@
 \global\x@false\global\xdimen@\z@
 \global\X@false\global\Xdimen@\z@
 \global\y@false\global\ydimen@\z@
 \global\Y@false\global\Ydimen@\z@
 \global\p@false\global\pdimen@\z@
 \global\label@ifalse\global\label@iifalse
 \global\dl@ifalse\global\ldimen@i\z@
 \global\dl@iifalse\global\ldimen@ii\z@
 \global\short@false\global\unshort@false}
\newif\iflabel@i
\newif\iflabel@ii
\newcount\scount@
\newcount\tcount@
\newcount\arrcount@
\newif\ifs@
\newdimen\sxdimen@
\newdimen\sydimen@
\newif\iftX@
\newdimen\tXdimen@i
\newdimen\tXdimen@ii
\newif\iftY@
\newdimen\tYdimen@i
\newdimen\tYdimen@ii
\newif\ifa@
\newcount\exacount@
\newif\ifx@
\newdimen\xdimen@
\newif\ifX@
\newdimen\Xdimen@
\newif\ify@
\newdimen\ydimen@
\newif\ifY@
\newdimen\Ydimen@
\newif\ifp@
\newdimen\pdimen@
\newif\ifdl@i
\newif\ifdl@ii
\newdimen\ldimen@i
\newdimen\ldimen@ii
\newif\ifshort@
\newif\ifunshort@
\def\zero@#1{\ifnum\scount@=\z@
 \if#1e\global\scount@\m@ne\else
 \if#1t\global\scount@\tw@\else
 \if#1h\global\scount@\thr@@\else
 \if#1'\global\scount@6 \else
 \if#1`\global\scount@7 \else
 \if#1(\global\scount@8 \else
 \if#1)\global\scount@9 \else
 \if#1s\global\scount@12 \else
 \if#1H\global\scount@13 \else
 \Err@{\Invalid@@ option \string\0}\fi\fi\fi\fi\fi\fi\fi\fi\fi
 \fi}
\def\one@#1{\ifnum\tcount@=\z@
 \if#1e\global\tcount@\m@ne\else
 \if#1h\global\tcount@\tw@\else
 \if#1t\global\tcount@\thr@@\else
 \if#1'\global\tcount@4 \else
 \if#1`\global\tcount@5 \else
 \if#1(\global\tcount@10 \else
 \if#1)\global\tcount@11 \else
 \if#1s\global\tcount@12 \else
 \if#1H\global\tcount@13 \else
 \Err@{\Invalid@@ option \string\1}\fi\fi\fi\fi\fi\fi\fi\fi\fi
 \fi}
\def\a@#1{\ifnum\arrcount@=\z@
 \if#10\global\arrcount@\m@ne\else
 \if#1+\global\arrcount@\@ne\else
 \if#1-\global\arrcount@\tw@\else
 \if#1=\global\arrcount@\thr@@\else
 \Err@{\Invalid@@ option \string\a}\fi\fi\fi\fi
 \fi}
\def\ds@(#1;#2){\ifs@\else
 \global\s@true
 \sxdimen@\hunit \global\sxdimen@#1\sxdimen@\relax
 \sydimen@\vunit \global\sydimen@#2\sydimen@\relax
 \fi}
\def\dtX@(#1;#2){\iftX@\else
 \global\tX@true
 \tXdimen@i\hunit \global\tXdimen@i#1\tXdimen@i\relax
 \tXdimen@ii\vunit \global\tXdimen@ii#2\tXdimen@ii\relax
 \fi}
\def\dtY@(#1;#2){\iftY@\else
 \global\tY@true
 \tYdimen@i\hunit \global\tYdimen@i#1\tYdimen@i\relax
 \tYdimen@ii\vunit \global\tYdimen@ii#2\tYdimen@ii\relax
 \fi}
\def\da@#1{\ifa@\else\global\a@true\global\exacount@#1\relax\fi}
\def\dx@#1{\ifx@\else
 \global\x@true
 \xdimen@\hunit \global\xdimen@#1\xdimen@\relax
 \fi}
\def\dX@#1{\ifX@\else
 \global\X@true
 \Xdimen@\hunit \global\Xdimen@#1\Xdimen@\relax
 \fi}
\def\dy@#1{\ify@\else
 \global\y@true
 \ydimen@\vunit \global\ydimen@#1\ydimen@\relax
 \fi}
\def\dY@#1{\ifY@\else
 \global\Y@true
 \Ydimen@\vunit \global\Ydimen@#1\Ydimen@\relax
 \fi}
\def\p@@#1{\ifp@\else
 \global\p@true
 \pdimen@\hunit \divide\pdimen@\tw@ \global\pdimen@#1\pdimen@\relax
 \fi}
\def\L@#1{\iflabel@i\else
 \global\label@itrue \gdef\label@i{#1}%
 \fi}
\def\l@#1{\iflabel@ii\else
 \global\label@iitrue \gdef\label@ii{#1}%
 \fi}
\def\dL@#1{\ifdl@i\else
 \global\dl@itrue \ldimen@i\hunit \global\ldimen@i#1\ldimen@i\relax
 \fi}
\def\dl@#1{\ifdl@ii\else
 \global\dl@iitrue \ldimen@ii\hunit \global\ldimen@ii#1\ldimen@ii\relax
 \fi}
\def\s@{\ifunshort@\else\global\short@true\fi}
\def\uns@{\ifshort@\else\global\unshort@true\global\short@false\fi}
\def\optioncodes@{\let\0\zero@\let\1\one@\let\a\a@\let\ds\ds@\let\dtX\dtX@
 \let\dtY\dtY@\let\da\da@\let\dx\dx@\let\dX\dX@\let\dY\dY@\let\dy\dy@
 \let\p\p@@\let\L\L@\let\l\l@\let\dL\dL@\let\dl\dl@\let\s\s@\let\uns\uns@}
\def\slopes@{\\161\\152\\143\\134\\255\\126\\357\\238\\349\\45{10}\\56{11}%
 \\11{12}\\65{13}\\54{14}\\43{15}\\32{16}\\53{17}\\21{18}\\52{19}\\31{20}%
 \\41{21}\\51{22}\\61{23}}
\newcount\tan@i
\newcount\tan@ip
\newcount\tan@ii
\newcount\tan@iip
\newdimen\slope@i
\newdimen\slope@ip
\newdimen\slope@ii
\newdimen\slope@iip
\newcount\angcount@
\newcount\extracount@
\def\slope@{{\slope@i=\secondy@ \advance\slope@i-\firsty@
 \ifN@\else\multiply\slope@i\m@ne\fi
 \slope@ii=\secondx@ \advance\slope@ii-\firstx@
 \ifE@\else\multiply\slope@ii\m@ne\fi
 \ifdim\slope@ii<\z@
  \global\tan@i6 \global\tan@ii\@ne \global\angcount@23
 \else
  \dimen@\slope@i \multiply\dimen@6
  \ifdim\dimen@<\slope@ii
   \global\tan@i\@ne \global\tan@ii6 \global\angcount@\@ne
  \else
   \dimen@\slope@ii \multiply\dimen@6
   \ifdim\dimen@<\slope@i
    \global\tan@i6 \global\tan@ii\@ne \global\angcount@23
   \else
    \tan@ip\z@ \tan@iip \@ne
    \def\\##1##2##3{\global\angcount@=##3\relax
     \slope@ip\slope@i \slope@iip\slope@ii
     \multiply\slope@iip##1\relax \multiply\slope@ip##2\relax
     \ifdim\slope@iip<\slope@ip
      \global\tan@ip=##1\relax \global\tan@iip=##2\relax
     \else
      \global\tan@i=##1\relax \global\tan@ii=##2\relax
      \def\\####1####2####3{}%
     \fi}%
    \slopes@
    \slope@i=\secondy@ \advance\slope@i-\firsty@
    \ifN@\else\multiply\slope@i\m@ne\fi
    \multiply\slope@i\tan@ii \multiply\slope@i\tan@iip \multiply\slope@i\tw@
    \count@\tan@i \multiply\count@\tan@iip
    \extracount@\tan@ip \multiply\extracount@\tan@ii
    \advance\count@\extracount@
    \slope@ii=\secondx@ \advance\slope@ii-\firstx@
    \ifE@\else\multiply\slope@ii\m@ne\fi
    \multiply\slope@ii\count@
    \ifdim\slope@i<\slope@ii
     \global\tan@i=\tan@ip \global\tan@ii=\tan@iip
     \global\advance\angcount@\m@ne
    \fi
   \fi
  \fi
 \fi}%
}
\def\slope@a#1{{\def\\##1##2##3{\ifnum##3=#1\global\tan@i=##1\relax
 \global\tan@ii=##2\relax\fi}\slopes@}}
\newcount\i@
\newcount\j@
\newcount\colcount@
\newcount\Colcount@
\newcount\tcolcount@
\newdimen\rowht@
\newdimen\rowdp@
\newcount\rowcount@
\newcount\Rowcount@
\newcount\maxcolrow@
\newtoks\colwidthtoks@
\newtoks\Rowheighttoks@
\newtoks\Rowdepthtoks@
\newtoks\widthtoks@
\newtoks\Widthtoks@
\newtoks\heighttoks@
\newtoks\Heighttoks@
\newtoks\depthtoks@
\newtoks\Depthtoks@
\newif\iffirstnewCDcr@
\def\dotoks@i{%
 \global\widthtoks@=\expandafter{\the\widthtoks@\else\getdim@\z@\fi}%
 \global\heighttoks@=\expandafter{\the\heighttoks@\else\getdim@\z@\fi}%
 \global\depthtoks@=\expandafter{\the\depthtoks@\else\getdim@\z@\fi}}
\def\dotoks@ii{%
 \global\widthtoks@{\ifcase\j@}%
 \global\heighttoks@{\ifcase\j@}%
 \global\depthtoks@{\ifcase\j@}}
\def\prenewCD@#1\endnewCD{\setbox\zer@
 \vbox{%
  \def\arrow@##1##2{{}}%
  \rowcount@\m@ne \colcount@\z@ \Colcount@\z@
  \firstnewCDcr@true \toks@{}%
  \widthtoks@{\ifcase\j@}%
  \Widthtoks@{\ifcase\i@}%
  \heighttoks@{\ifcase\j@}%
  \Heighttoks@{\ifcase\i@}%
  \depthtoks@{\ifcase\j@}%
  \Depthtoks@{\ifcase\i@}%
  \Rowheighttoks@{\ifcase\i@}%
  \Rowdepthtoks@{\ifcase\i@}%
  \Let@
  \everycr{%
   \noalign{%
    \global\advance\rowcount@\@ne
    \ifnum\colcount@<\Colcount@
    \else
     \global\Colcount@=\colcount@ \global\maxcolrow@=\rowcount@
    \fi
    \global\colcount@\z@
    \iffirstnewCDcr@
     \global\firstnewCDcr@false
    \else
     \edef\next@{\the\Rowheighttoks@\noexpand\or\noexpand\getdim@\the\rowht@}%
      \global\Rowheighttoks@=\expandafter{\next@}%
     \edef\next@{\the\Rowdepthtoks@\noexpand\or\noexpand\getdim@\the\rowdp@}%
      \global\Rowdepthtoks@=\expandafter{\next@}%
     \global\rowht@\z@ \global\rowdp@\z@
     \dotoks@i
     \edef\next@{\the\Widthtoks@\noexpand\or\the\widthtoks@}%
      \global\Widthtoks@=\expandafter{\next@}%
     \edef\next@{\the\Heighttoks@\noexpand\or\the\heighttoks@}%
      \global\Heighttoks@=\expandafter{\next@}%
     \edef\next@{\the\Depthtoks@\noexpand\or\the\depthtoks@}%
      \global\Depthtoks@=\expandafter{\next@}%
     \dotoks@ii
    \fi}}%
  \tabskip\z@
  \halign{&\setbox\zer@\hbox{\vrule height10\p@ width\z@ depth\z@
   $\m@th\displaystyle{##}$}\copy\zer@
   \ifdim\ht\zer@>\rowht@ \global\rowht@\ht\zer@ \fi
   \ifdim\dp\zer@>\rowdp@ \global\rowdp@\dp\zer@ \fi
   \global\advance\colcount@\@ne
   \edef\next@{\the\widthtoks@\noexpand\or\noexpand\getdim@\the\wd\zer@}%
    \global\widthtoks@=\expandafter{\next@}%
   \edef\next@{\the\heighttoks@\noexpand\or\noexpand\getdim@\the\ht\zer@}%
    \global\heighttoks@=\expandafter{\next@}%
   \edef\next@{\the\depthtoks@\noexpand\or\noexpand\getdim@\the\dp\zer@}%
    \global\depthtoks@=\expandafter{\next@}%
   \cr#1\crcr}}%
 \Rowcount@=\rowcount@
 \global\Widthtoks@=\expandafter{\the\Widthtoks@\fi\relax}%
 \edef\Width@##1##2{\i@=##1\relax\j@=##2\relax\the\Widthtoks@}%
 \global\Heighttoks@=\expandafter{\the\Heighttoks@\fi\relax}%
 \edef\Height@##1##2{\i@=##1\relax\j@=##2\relax\the\Heighttoks@}%
 \global\Depthtoks@=\expandafter{\the\Depthtoks@\fi\relax}%
 \edef\Depth@##1##2{\i@=##1\relax\j@=##2\relax\the\Depthtoks@}%
 \edef\next@{\the\Rowheighttoks@\noexpand\fi\relax}%
 \global\Rowheighttoks@=\expandafter{\next@}%
 \edef\Rowheight@##1{\i@=##1\relax\the\Rowheighttoks@}%
 \edef\next@{\the\Rowdepthtoks@\noexpand\fi\relax}%
 \global\Rowdepthtoks@=\expandafter{\next@}%
 \edef\Rowdepth@##1{\i@=##1\relax\the\Rowdepthtoks@}%
 \colwidthtoks@{\fi}%
 \setbox\zer@\vbox{%
  \unvbox\zer@
  \count@\rowcount@
  \loop
   \unskip\unpenalty
   \setbox\zer@\lastbox
   \ifnum\count@>\maxcolrow@ \advance\count@\m@ne
   \repeat
  \hbox{%
   \unhbox\zer@
   \count@\z@
   \loop
    \unskip
    \setbox\zer@\lastbox
    \edef\next@{\noexpand\or\noexpand\getdim@\the\wd\zer@\the\colwidthtoks@}%
     \global\colwidthtoks@=\expandafter{\next@}%
    \advance\count@\@ne
    \ifnum\count@<\Colcount@
    \repeat}}%
 \edef\next@{\noexpand\ifcase\noexpand\i@\the\colwidthtoks@}%
  \global\colwidthtoks@=\expandafter{\next@}%
 \edef\Colwidth@##1{\i@=##1\relax\the\colwidthtoks@}%
 \colwidthtoks@{}\Rowheighttoks@{}\Rowdepthtoks@{}\widthtoks@{}%
 \Widthtoks@{}\heighttoks@{}\Heighttoks@{}\depthtoks@{}\Depthtoks@{}%
}
\newcount\xoff@
\newcount\yoff@
\newcount\endcount@
\newcount\rcount@
\newdimen\firstx@
\newdimen\firsty@
\newdimen\secondx@
\newdimen\secondy@
\newdimen\tocenter@
\newdimen\charht@
\newdimen\charwd@
\def\outside@{\Err@{This arrow points outside the \string\newCD}}
\newif\ifsvertex@
\newif\iftvertex@
\def\arrow@#1#2{\xoff@=#1\relax\yoff@=#2\relax
 \count@\rowcount@ \advance\count@-\yoff@
 \ifnum\count@<\@ne \outside@ \else \ifnum\count@>\Rowcount@ \outside@ \fi\fi
 \count@\colcount@ \advance\count@\xoff@
 \ifnum\count@<\@ne \outside@ \else \ifnum\count@>\Colcount@ \outside@\fi\fi
 \tcolcount@\colcount@ \advance\tcolcount@\xoff@
 \Width@\rowcount@\colcount@ \tocenter@=-\getdim@ \divide\tocenter@\tw@
 \ifdim\getdim@=\z@
  \firstx@\z@ \firsty@\mathaxis@ \svertex@true
 \else
  \svertex@false
  \ifHshort@
   \Colwidth@\colcount@
    \ifE@ \firstx@=.5\getdim@ \else \firstx@=-.5\getdim@ \fi
  \else
   \ifE@ \firstx@=\getdim@ \else \firstx@=-\getdim@ \fi
   \divide\firstx@\tw@
  \fi
  \ifE@
   \ifH@ \advance\firstx@\thr@@\p@ \else \advance\firstx@-\thr@@\p@ \fi
  \else
   \ifH@ \advance\firstx@-\thr@@\p@ \else \advance\firstx@\thr@@\p@ \fi
  \fi
  \ifN@
   \Height@\rowcount@\colcount@ \firsty@=\getdim@
   \ifV@ \advance\firsty@\thr@@\p@ \fi
  \else
   \ifV@
    \Depth@\rowcount@\colcount@ \firsty@=-\getdim@
    \advance\firsty@-\thr@@\p@
   \else
    \firsty@\z@
   \fi
  \fi
 \fi
 \ifV@
 \else
  \Colwidth@\colcount@
  \ifE@ \secondx@=\getdim@ \else \secondx@=-\getdim@ \fi
  \divide\secondx@\tw@
  \ifE@ \else \getcgap@\colcount@ \advance\secondx@-\getdim@ \fi
  \endcount@=\colcount@ \advance\endcount@\xoff@
  \count@=\colcount@
  \ifE@
   \advance\count@\@ne
   \loop
    \ifnum\count@<\endcount@
    \Colwidth@\count@ \advance\secondx@\getdim@
    \getcgap@\count@ \advance\secondx@\getdim@
    \advance\count@\@ne
    \repeat
  \else
   \advance\count@\m@ne
   \loop
    \ifnum\count@>\endcount@
    \Colwidth@\count@ \advance\secondx@-\getdim@
    \getcgap@\count@ \advance\secondx@-\getdim@
    \advance\count@\m@ne
    \repeat
  \fi
  \Colwidth@\count@ \divide\getdim@\tw@
  \ifHshort@
  \else
   \ifE@ \advance\secondx@\getdim@ \else \advance\secondx@-\getdim@ \fi
  \fi
  \ifE@ \getcgap@\count@ \advance\secondx@\getdim@ \fi
  \rcount@\rowcount@ \advance\rcount@-\yoff@
  \Width@\rcount@\count@ \divide\getdim@\tw@
  \tvertex@false
  \ifH@\ifdim\getdim@=\z@\tvertex@true\Hshort@false\fi\fi
  \ifHshort@
  \else
   \ifE@ \advance\secondx@-\getdim@ \else \advance\secondx@\getdim@ \fi
  \fi
  \iftvertex@
   \advance\secondx@.4\p@
  \else
   \ifE@ \advance\secondx@-\thr@@\p@ \else \advance\secondx@\thr@@\p@ \fi
  \fi
 \fi
 \ifH@
 \else
  \ifN@
   \Rowheight@\rowcount@ \secondy@\getdim@
  \else
   \Rowdepth@\rowcount@ \secondy@-\getdim@
   \getrgap@\rowcount@ \advance\secondy@-\getdim@
  \fi
  \endcount@=\rowcount@ \advance\endcount@-\yoff@
  \count@=\rowcount@
  \ifN@
   \advance\count@\m@ne
   \loop
    \ifnum\count@>\endcount@
    \Rowheight@\count@ \advance\secondy@\getdim@
    \Rowdepth@\count@ \advance\secondy@\getdim@
    \getrgap@\count@ \advance\secondy@\getdim@
    \advance\count@\m@ne
    \repeat
  \else
   \advance\count@\@ne
   \loop
    \ifnum\count@<\endcount@
    \Rowheight@\count@ \advance\secondy@-\getdim@
    \Rowdepth@\count@ \advance\secondy@-\getdim@
    \getrgap@\count@ \advance\secondy@-\getdim@
    \advance\count@\@ne
    \repeat
  \fi
  \tvertex@false
  \ifV@\Width@\count@\colcount@\ifdim\getdim@=\z@\tvertex@true\fi\fi
  \ifN@
   \getrgap@\count@ \advance\secondy@\getdim@
   \Rowdepth@\count@ \advance\secondy@\getdim@
   \iftvertex@
    \advance\secondy@\mathaxis@
   \else
    \Depth@\count@\tcolcount@ \advance\secondy@-\getdim@
    \advance\secondy@-\thr@@\p@
   \fi
  \else
   \Rowheight@\count@ \advance\secondy@-\getdim@
   \iftvertex@
    \advance\secondy@\mathaxis@
   \else
    \Height@\count@\tcolcount@ \advance\secondy@\getdim@
    \advance\secondy@\thr@@\p@
   \fi
  \fi
 \fi
 \ifV@\else\advance\firstx@\sxdimen@\fi
 \ifH@\else\advance\firsty@\sydimen@\fi
 \iftX@
  \advance\secondy@\tXdimen@ii
  \advance\secondx@\tXdimen@i
  \slope@
 \else
  \iftY@
   \advance\secondy@\tYdimen@ii
   \advance\secondx@\tYdimen@i
   \slope@
   \secondy@=\secondx@ \advance\secondy@-\firstx@
   \ifNESW@ \else \multiply\secondy@\m@ne \fi
   \multiply\secondy@\tan@i \divide\secondy@\tan@ii \advance\secondy@\firsty@
  \else
   \ifa@
    \slope@
    \ifNESW@ \global\advance\angcount@\exacount@ \else
      \global\advance\angcount@-\exacount@ \fi
    \ifnum\angcount@>23 \angcount@23 \fi
    \ifnum\angcount@<\@ne \angcount@\@ne \fi
    \slope@a\angcount@
    \ifY@
     \advance\secondy@\Ydimen@
    \else
     \ifX@
      \advance\secondx@\Xdimen@
      \dimen@\secondx@ \advance\dimen@-\firstx@
      \ifNESW@\else\multiply\dimen@\m@ne\fi
      \multiply\dimen@\tan@i \divide\dimen@\tan@ii
      \advance\dimen@\firsty@ \secondy@=\dimen@
     \fi
    \fi
   \else
    \ifH@\else\ifV@\else\slope@\fi\fi
   \fi
  \fi
 \fi
 \ifH@\else\ifV@\else\ifsvertex@\else
  \dimen@=6\p@ \multiply\dimen@\tan@ii
  \count@=\tan@i \advance\count@\tan@ii \divide\dimen@\count@
  \ifE@ \advance\firstx@\dimen@ \else \advance\firstx@-\dimen@ \fi
  \multiply\dimen@\tan@i \divide\dimen@\tan@ii
  \ifN@ \advance\firsty@\dimen@ \else \advance\firsty@-\dimen@ \fi
 \fi\fi\fi
 \ifp@
  \ifH@\else\ifV@\else
   \getcos@\pdimen@ \advance\firsty@\dimen@ \advance\secondy@\dimen@
   \ifNESW@ \advance\firstx@-\dimen@ii \else \advance\firstx@\dimen@ii \fi
  \fi\fi
 \fi
 \ifH@\else\ifV@\else
  \ifnum\tan@i>\tan@ii
   \charht@=10\p@ \charwd@=10\p@
   \multiply\charwd@\tan@ii \divide\charwd@\tan@i
  \else
   \charwd@=10\p@ \charht@=10\p@
   \divide\charht@\tan@ii \multiply\charht@\tan@i
  \fi
  \ifnum\tcount@=\thr@@
   \ifN@ \advance\secondy@-.3\charht@ \else\advance\secondy@.3\charht@ \fi
  \fi
  \ifnum\scount@=\tw@
   \ifE@ \advance\firstx@.3\charht@ \else \advance\firstx@-.3\charht@ \fi
  \fi
  \ifnum\tcount@=12
   \ifN@ \advance\secondy@-\charht@ \else \advance\secondy@\charht@ \fi
  \fi
  \iftY@
  \else
   \ifa@
    \ifX@
    \else
     \secondx@\secondy@ \advance\secondx@-\firsty@
     \ifNESW@\else\multiply\secondx@\m@ne\fi
     \multiply\secondx@\tan@ii \divide\secondx@\tan@i
     \advance\secondx@\firstx@
    \fi
   \fi
  \fi
 \fi\fi
 \ifH@\harrow@\else\ifV@\varrow@\else\arrow@@\fi\fi}
\newdimen\mathaxis@
\mathaxis@90\p@ \divide\mathaxis@36
\def\harrow@b{\ifE@\hskip\tocenter@\hskip\firstx@\fi}
\def\harrow@bb{\ifE@\hskip\xdimen@\else\hskip\Xdimen@\fi}
\def\harrow@e{\ifE@\else\hskip-\firstx@\hskip-\tocenter@\fi}
\def\harrow@ee{\ifE@\hskip-\Xdimen@\else\hskip-\xdimen@\fi}
\def\harrow@{\dimen@\secondx@\advance\dimen@-\firstx@
 \ifE@ \let\next@\rlap \else  \multiply\dimen@\m@ne \let\next@\llap \fi
 \next@{%
  \harrow@b
  \smash{\raise\pdimen@\hbox to\dimen@
   {\harrow@bb\arrow@ii
    \ifnum\arrcount@=\m@ne \else \ifnum\arrcount@=\thr@@ \else
     \ifE@
      \ifnum\scount@=\m@ne
      \else
       \ifcase\scount@\or\or\char118 \or\char117 \or\or\or\char119 \or
       \char120 \or\char121 \or\char122 \or\or\or\arrow@i\char125 \or
       \char117 \hskip\thr@@\p@\char117 \hskip-\thr@@\p@\fi
      \fi
     \else
      \ifnum\tcount@=\m@ne
      \else
       \ifcase\tcount@\char117 \or\or\char117 \or\char118 \or\char119 \or
       \char120\or\or\or\or\or\char121 \or\char122 \or\arrow@i\char125
       \or\char117 \hskip\thr@@\p@\char117 \hskip-\thr@@\p@\fi
      \fi
     \fi
    \fi\fi
    \dimen@\mathaxis@ \advance\dimen@.2\p@
    \dimen@ii\mathaxis@ \advance\dimen@ii-.2\p@
    \ifnum\arrcount@=\m@ne
     \let\leads@\null
    \else
     \ifcase\arrcount@
      \def\leads@{\hrule height\dimen@ depth-\dimen@ii}\or
      \def\leads@{\hrule height\dimen@ depth-\dimen@ii}\or
      \def\leads@{\hbox to10\p@{%
       \leaders\hrule height\dimen@ depth-\dimen@ii\hfil
       \hfil
      \leaders\hrule height\dimen@ depth-\dimen@ii\hskip\z@ plus2fil\relax
       \hfil
       \leaders\hrule height\dimen@ depth-\dimen@ii\hfil}}\or
     \def\leads@{\hbox{\hbox to10\p@{\dimen@\mathaxis@ \advance\dimen@1.2\p@
       \dimen@ii\dimen@ \advance\dimen@ii-.4\p@
       \leaders\hrule height\dimen@ depth-\dimen@ii\hfil}%
       \kern-10\p@
       \hbox to10\p@{\dimen@\mathaxis@ \advance\dimen@-1.2\p@
       \dimen@ii\dimen@ \advance\dimen@ii-.4\p@
       \leaders\hrule height\dimen@ depth-\dimen@ii\hfil}}}\fi
    \fi
    \cleaders\leads@\hfil
    \ifnum\arrcount@=\m@ne\else\ifnum\arrcount@=\thr@@\else
     \arrow@i
     \ifE@
      \ifnum\tcount@=\m@ne
      \else
       \ifcase\tcount@\char119 \or\or\char119 \or\char120 \or\char121 \or
       \char122 \or \or\or\or\or\char123\or\char124 \or
       \char125 \or\char119 \hskip-\thr@@\p@\char119 \hskip\thr@@\p@\fi
      \fi
     \else
      \ifcase\scount@\or\or\char120 \or\char119 \or\or\or\char121 \or\char122
      \or\char123 \or\char124 \or\or\or\char125 \or
      \char119 \hskip-\thr@@\p@\char119 \hskip\thr@@\p@\fi
     \fi
    \fi\fi
    \harrow@ee}}%
  \harrow@e}%
 \iflabel@i
  \dimen@ii\z@ \setbox\zer@\hbox{$\m@th\tsize@@\label@i$}%
  \ifnum\arrcount@=\m@ne
  \else
   \advance\dimen@ii\mathaxis@
   \advance\dimen@ii\dp\zer@ \advance\dimen@ii\tw@\p@
   \ifnum\arrcount@=\thr@@ \advance\dimen@ii\tw@\p@ \fi
  \fi
  \advance\dimen@ii\pdimen@
  \next@{\harrow@b\smash{\raise\dimen@ii\hbox to\dimen@
   {\harrow@bb\hskip\tw@\ldimen@i\hfil\box\zer@\hfil\harrow@ee}}\harrow@e}%
 \fi
 \iflabel@ii
  \ifnum\arrcount@=\m@ne
  \else
   \setbox\zer@\hbox{$\m@th\tsize@\label@ii$}%
   \dimen@ii-\ht\zer@ \advance\dimen@ii-\tw@\p@
   \ifnum\arrcount@=\thr@@ \advance\dimen@ii-\tw@\p@ \fi
   \advance\dimen@ii\mathaxis@ \advance\dimen@ii\pdimen@
   \next@{\harrow@b\smash{\raise\dimen@ii\hbox to\dimen@
    {\harrow@bb\hskip\tw@\ldimen@ii\hfil\box\zer@\hfil\harrow@ee}}\harrow@e}%
  \fi
 \fi}
\let\tsize@\tsize
\def\tsizenewCDlabels{\let\tsize@\tsize}
\def\ssizenewCDlabels{\let\tsize@\ssize}
\def\tsize@@{\ifnum\arrcount@=\m@ne\else\tsize@\fi}
\def\varrow@{\dimen@\secondy@ \advance\dimen@-\firsty@
 \ifN@ \else \multiply\dimen@\m@ne \fi
 \setbox\zer@\vbox to\dimen@
  {\ifN@ \vskip-\Ydimen@ \else \vskip\ydimen@ \fi
   \ifnum\arrcount@=\m@ne\else\ifnum\arrcount@=\thr@@\else
    \hbox{\arrow@iii
     \ifN@
      \ifnum\tcount@=\m@ne
      \else
       \ifcase\tcount@\char117 \or\or\char117 \or\char118 \or\char119 \or
       \char120 \or\or\or\or\or\char121 \or\char122 \or\char123 \or
       \vbox{\hbox{\char117 }\nointerlineskip\vskip\thr@@\p@
       \hbox{\char117 }\vskip-\thr@@\p@}\fi
      \fi
     \else
      \ifcase\scount@\or\or\char118 \or\char117 \or\or\or\char119 \or
      \char120 \or\char121 \or\char122 \or\or\or\char123 \or
      \vbox{\hbox{\char117 }\nointerlineskip\vskip\thr@@\p@
      \hbox{\char117 }\vskip-\thr@@\p@}\fi
     \fi}%
    \nointerlineskip
   \fi\fi
   \ifnum\arrcount@=\m@ne
    \let\leads@\null
   \else
    \ifcase\arrcount@\let\leads@\vrule\or\let\leads@\vrule\or
    \def\leads@{\vbox to10\p@{%
     \hrule height 1.67\p@ depth\z@ width.4\p@
     \vfil
     \hrule height 3.33\p@ depth\z@ width.4\p@
     \vfil
     \hrule height 1.67\p@ depth\z@ width.4\p@}}\or
    \def\leads@{\hbox{\vrule height\p@\hskip\tw@\p@\vrule}}\fi
   \fi
  \cleaders\leads@\vfill\nointerlineskip
   \ifnum\arrcount@=\m@ne\else\ifnum\arrcount@=\thr@@\else
    \hbox{\arrow@iv
     \ifN@
      \ifcase\scount@\or\or\char118 \or\char117 \or\or\or\char119 \or
      \char120 \or\char121 \or\char122 \or\or\or\arrow@iii\char123 \or
      \vbox{\hbox{\char117 }\nointerlineskip\vskip-\thr@@\p@
      \hbox{\char117 }\vskip\thr@@\p@}\fi
     \else
      \ifnum\tcount@=\m@ne
      \else
       \ifcase\tcount@\char117 \or\or\char117 \or\char118 \or\char119 \or
       \char120 \or\or\or\or\or\char121 \or\char122 \or\arrow@iii\char123 \or
       \vbox{\hbox{\char117 }\nointerlineskip\vskip-\thr@@\p@
       \hbox{\char117 }\vskip\thr@@\p@}\fi
      \fi
     \fi}%
   \fi\fi
   \ifN@\vskip\ydimen@\else\vskip-\Ydimen@\fi}%
 \ifN@
  \dimen@ii\firsty@
 \else
  \dimen@ii-\firsty@ \advance\dimen@ii\ht\zer@ \multiply\dimen@ii\m@ne
 \fi
 \rlap{\smash{\hskip\tocenter@ \hskip\pdimen@ \raise\dimen@ii \box\zer@}}%
 \iflabel@i
  \setbox\zer@\vbox to\dimen@{\vfil
   \hbox{$\m@th\tsize@@\label@i$}\vskip\tw@\ldimen@i\vfil}%
  \rlap{\smash{\hskip\tocenter@ \hskip\pdimen@
  \ifnum\arrcount@=\m@ne \let\next@\relax \else \let\next@\llap \fi
  \next@{\raise\dimen@ii\hbox{\ifnum\arrcount@=\m@ne \hskip-.5\wd\zer@ \fi
   \box\zer@ \ifnum\arrcount@=\m@ne \else \hskip\tw@\p@ \fi}}}}%
 \fi
 \iflabel@ii
  \ifnum\arrcount@=\m@ne
  \else
   \setbox\zer@\vbox to\dimen@{\vfil
    \hbox{$\m@th\tsize@\label@ii$}\vskip\tw@\ldimen@ii\vfil}%
   \rlap{\smash{\hskip\tocenter@ \hskip\pdimen@
   \rlap{\raise\dimen@ii\hbox{\ifnum\arrcount@=\thr@@ \hskip4.5\p@ \else
    \hskip2.5\p@ \fi\box\zer@}}}}%
  \fi
 \fi
}
\newdimen\goal@
\newdimen\shifted@
\newcount\Tcount@
\newcount\Scount@
\newbox\shaft@
\newcount\slcount@
\def\getcos@#1{%
 \ifnum\tan@i<\tan@ii
  \dimen@#1%
  \ifnum\slcount@<8 \count@9 \else \ifnum\slcount@<12 \count@8 \else
   \count@7 \fi\fi
  \multiply\dimen@\count@ \divide\dimen@10
  \dimen@ii\dimen@ \multiply\dimen@ii\tan@i \divide\dimen@ii\tan@ii
 \else
  \dimen@ii#1%
  \count@-\slcount@ \advance\count@24
  \ifnum\count@<8 \count@9 \else \ifnum\count@<12 \count@8
   \else\count@7 \fi\fi
  \multiply\dimen@ii\count@ \divide\dimen@ii10
  \dimen@\dimen@ii \multiply\dimen@\tan@ii \divide\dimen@\tan@i
 \fi}
\newdimen\adjust@
\def\Nnext@{\ifN@\let\next@\raise\else\let\next@\lower\fi}
\def\arrow@@{\slcount@\angcount@
 \ifNESW@
  \ifnum\angcount@<10
   \let\arrowfont@=\arrow@i \advance\angcount@\m@ne \multiply\angcount@13
  \else
   \ifnum\angcount@<19
    \let\arrowfont@=\arrow@ii \advance\angcount@-10 \multiply\angcount@13
   \else
    \let\arrowfont@=\arrow@iii \advance\angcount@-19 \multiply\angcount@13
  \fi\fi
  \Tcount@\angcount@
 \else
  \ifnum\angcount@<5
   \let\arrowfont@=\arrow@iii \advance\angcount@\m@ne \multiply\angcount@13
   \advance\angcount@65
  \else
   \ifnum\angcount@<14
    \let\arrowfont@=\arrow@iv \advance\angcount@-5 \multiply\angcount@13
   \else
    \ifnum\angcount@<23
     \let\arrowfont@=\arrow@v \advance\angcount@-14 \multiply\angcount@13
    \else
     \let\arrowfont@=\arrow@i \angcount@=117
  \fi\fi\fi
  \ifnum\angcount@=117 \Tcount@=115 \else\Tcount@\angcount@ \fi
 \fi
 \Scount@\Tcount@
 \ifE@
  \ifnum\tcount@=\z@ \advance\Tcount@\tw@ \else\ifnum\tcount@=13
   \advance\Tcount@\tw@ \else \advance\Tcount@\tcount@ \fi\fi
  \ifnum\scount@=\z@ \else \ifnum\scount@=13 \advance\Scount@\thr@@ \else
   \advance\Scount@\scount@ \fi\fi
 \else
  \ifcase\tcount@\advance\Tcount@\thr@@\or\or\advance\Tcount@\thr@@\or
  \advance\Tcount@\tw@\or\advance\Tcount@6 \or\advance\Tcount@7
  \or\or\or\or\or \advance\Tcount@8 \or\advance\Tcount@9 \or
  \advance\Tcount@12 \or\advance\Tcount@\thr@@\fi
  \ifcase\scount@\or\or\advance\Scount@\thr@@\or\advance\Scount@\tw@\or
  \or\or\advance\Scount@4 \or\advance\Scount@5 \or\advance\Scount@10
  \or\advance\Scount@11 \or\or\or\advance\Scount@12 \or\advance
  \Scount@\tw@\fi
 \fi
 \ifcase\arrcount@\or\or\advance\angcount@\@ne\else\fi
 \ifN@ \shifted@=\firsty@ \else\shifted@=-\firsty@ \fi
 \ifE@ \else\advance\shifted@\charht@ \fi
 \goal@=\secondy@ \advance\goal@-\firsty@
 \ifN@\else\multiply\goal@\m@ne\fi
 \setbox\shaft@\hbox{\arrowfont@\char\angcount@}%
 \ifnum\arrcount@=\thr@@
  \getcos@{1.5\p@}%
  \setbox\shaft@\hbox to\wd\shaft@{\arrowfont@
   \rlap{\hskip\dimen@ii
    \smash{\ifNESW@\let\next@\lower\else\let\next@\raise\fi
     \next@\dimen@\hbox{\arrowfont@\char\angcount@}}}%
   \rlap{\hskip-\dimen@ii
    \smash{\ifNESW@\let\next@\raise\else\let\next@\lower\fi
      \next@\dimen@\hbox{\arrowfont@\char\angcount@}}}\hfil}%
 \fi
 \rlap{\smash{\hskip\tocenter@\hskip\firstx@
  \ifnum\arrcount@=\m@ne
  \else
   \ifnum\arrcount@=\thr@@
   \else
    \ifnum\scount@=\m@ne
    \else
     \ifnum\scount@=\z@
     \else
      \setbox\zer@\hbox{\ifnum\angcount@=117 \arrow@v\else\arrowfont@\fi
       \char\Scount@}%
      \ifNESW@
       \ifnum\scount@=\tw@
        \dimen@=\shifted@ \advance\dimen@-\charht@
        \ifN@\hskip-\wd\zer@\fi
        \Nnext@
        \next@\dimen@\copy\zer@
        \ifN@\else\hskip-\wd\zer@\fi
       \else
        \Nnext@
        \ifN@\else\hskip-\wd\zer@\fi
        \next@\shifted@\copy\zer@
        \ifN@\hskip-\wd\zer@\fi
       \fi
       \ifnum\scount@=12
        \advance\shifted@\charht@ \advance\goal@-\charht@
        \ifN@ \hskip\wd\zer@ \else \hskip-\wd\zer@ \fi
       \fi
       \ifnum\scount@=13
        \getcos@{\thr@@\p@}%
        \ifN@ \hskip\dimen@ \else \hskip-\wd\zer@ \hskip-\dimen@ \fi
        \adjust@\shifted@ \advance\adjust@\dimen@ii
        \Nnext@
        \next@\adjust@\copy\zer@
        \ifN@ \hskip-\dimen@ \hskip-\wd\zer@ \else \hskip\dimen@ \fi
       \fi
      \else
       \ifN@\hskip-\wd\zer@\fi
       \ifnum\scount@=\tw@
        \ifN@ \hskip\wd\zer@ \else \hskip-\wd\zer@ \fi
        \dimen@=\shifted@ \advance\dimen@-\charht@
        \Nnext@
        \next@\dimen@\copy\zer@
        \ifN@\hskip-\wd\zer@\fi
       \else
        \Nnext@
        \next@\shifted@\copy\zer@
        \ifN@\else\hskip-\wd\zer@\fi
       \fi
       \ifnum\scount@=12
        \advance\shifted@\charht@ \advance\goal@-\charht@
        \ifN@ \hskip-\wd\zer@ \else \hskip\wd\zer@ \fi
       \fi
       \ifnum\scount@=13
        \getcos@{\thr@@\p@}%
        \ifN@ \hskip-\wd\zer@ \hskip-\dimen@ \else \hskip\dimen@ \fi
        \adjust@\shifted@ \advance\adjust@\dimen@ii
        \Nnext@
        \next@\adjust@\copy\zer@
        \ifN@ \hskip\dimen@ \else \hskip-\dimen@ \hskip-\wd\zer@ \fi
       \fi	
      \fi
  \fi\fi\fi\fi
  \ifnum\arrcount@=\m@ne
  \else
   \loop
    \ifdim\goal@>\charht@
    \ifE@\else\hskip-\charwd@\fi
    \Nnext@
    \next@\shifted@\copy\shaft@
    \ifE@\else\hskip-\charwd@\fi
    \advance\shifted@\charht@ \advance\goal@ -\charht@
    \repeat
   \ifdim\goal@>\z@
    \dimen@=\charht@ \advance\dimen@-\goal@
    \divide\dimen@\tan@i \multiply\dimen@\tan@ii
    \ifE@ \hskip-\dimen@ \else \hskip-\charwd@ \hskip\dimen@ \fi
    \adjust@=\shifted@ \advance\adjust@-\charht@ \advance\adjust@\goal@
    \Nnext@
    \next@\adjust@\copy\shaft@
    \ifE@ \else \hskip-\charwd@ \fi
   \else
    \adjust@=\shifted@ \advance\adjust@-\charht@
   \fi
  \fi
  \ifnum\arrcount@=\m@ne
  \else
   \ifnum\arrcount@=\thr@@
   \else
    \ifnum\tcount@=\m@ne
    \else
     \setbox\zer@
      \hbox{\ifnum\angcount@=117 \arrow@v\else\arrowfont@\fi\char\Tcount@}%
     \ifnum\tcount@=\thr@@
      \advance\adjust@\charht@
      \ifE@\else\ifN@\hskip-\charwd@\else\hskip-\wd\zer@\fi\fi
     \else
      \ifnum\tcount@=12
       \advance\adjust@\charht@
       \ifE@\else\ifN@\hskip-\charwd@\else\hskip-\wd\zer@\fi\fi
      \else
       \ifE@\hskip-\wd\zer@\fi
     \fi\fi
     \Nnext@
     \next@\adjust@\copy\zer@
     \ifnum\tcount@=13
      \hskip-\wd\zer@
      \getcos@{\thr@@\p@}%
      \ifE@\hskip-\dimen@ \else\hskip\dimen@ \fi
      \advance\adjust@-\dimen@ii
      \Nnext@
      \next@\adjust@\box\zer@
     \fi
  \fi\fi\fi}}%
 \iflabel@i
  \rlap{\hskip\tocenter@
  \dimen@\firstx@ \advance\dimen@\secondx@ \divide\dimen@\tw@
  \advance\dimen@\ldimen@i
  \dimen@ii\firsty@ \advance\dimen@ii\secondy@ \divide\dimen@ii\tw@
  \multiply\ldimen@i\tan@i \divide\ldimen@i\tan@ii
  \ifNESW@ \advance\dimen@ii\ldimen@i \else \advance\dimen@ii-\ldimen@i \fi
  \setbox\zer@\hbox{\ifNESW@\else\ifnum\arrcount@=\thr@@\hskip4\p@\else
   \hskip\tw@\p@\fi\fi
   $\m@th\tsize@@\label@i$\ifNESW@\ifnum\arrcount@=\thr@@\hskip4\p@\else
   \hskip\tw@\p@\fi\fi}%
  \ifnum\arrcount@=\m@ne
   \ifNESW@ \advance\dimen@.5\wd\zer@ \advance\dimen@\p@ \else
    \advance\dimen@-.5\wd\zer@ \advance\dimen@-\p@ \fi
   \advance\dimen@ii-.5\ht\zer@
  \else
   \advance\dimen@ii\dp\zer@
   \ifnum\slcount@<6 \advance\dimen@ii\tw@\p@ \fi
  \fi
  \hskip\dimen@
  \ifNESW@ \let\next@\llap \else\let\next@\rlap \fi
  \next@{\smash{\raise\dimen@ii\box\zer@}}}%
 \fi
 \iflabel@ii
  \ifnum\arrcount@=\m@ne
  \else
   \rlap{\hskip\tocenter@
   \dimen@\firstx@ \advance\dimen@\secondx@ \divide\dimen@\tw@
   \ifNESW@ \advance\dimen@\ldimen@ii \else \advance\dimen@-\ldimen@ii \fi
   \dimen@ii\firsty@ \advance\dimen@ii\secondy@ \divide\dimen@ii\tw@
   \multiply\ldimen@ii\tan@i \divide\ldimen@ii\tan@ii
   \advance\dimen@ii\ldimen@ii
   \setbox\zer@\hbox{\ifNESW@\ifnum\arrcount@=\thr@@\hskip4\p@\else
    \hskip\tw@\p@\fi\fi
    $\m@th\tsize@\label@ii$\ifNESW@\else\ifnum\arrcount@=\thr@@\hskip4\p@
    \else\hskip\tw@\p@\fi\fi}%
   \advance\dimen@ii-\ht\zer@
   \ifnum\slcount@<9 \advance\dimen@ii-\thr@@\p@ \fi
   \ifNESW@ \let\next@\rlap \else \let\next@\llap \fi
   \hskip\dimen@\next@{\smash{\raise\dimen@ii\box\zer@}}}%
  \fi
 \fi
}
\def\outnewCD@#1{\def#1{\Err@{\string#1 must not be used within \string\newCD}}}
\newskip\prenewCDskip@
\newskip\postnewCDskip@
\prenewCDskip@\z@
\postnewCDskip@\z@
\def\prenewCDspace#1{\RIfMIfI@
 \onlydmatherr@\prenewCDspace\else\advance\prenewCDskip@#1\relax\fi\else
 \onlydmatherr@\prenewCDspace\fi}
\def\postnewCDspace#1{\RIfMIfI@
 \onlydmatherr@\postnewCDspace\else\advance\postnewCDskip@#1\relax\fi\else
 \onlydmatherr@\postnewCDspace\fi}
\def\predisplayspace#1{\RIfMIfI@
 \onlydmatherr@\predisplayspace\else
 \advance\abovedisplayskip#1\relax
 \advance\abovedisplayshortskip#1\relax\fi
 \else\onlydmatherr@\prenewCDspace\fi}
\def\postdisplayspace#1{\RIfMIfI@
 \onlydmatherr@\postdisplayspace\else
 \advance\belowdisplayskip#1\relax
 \advance\belowdisplayshortskip#1\relax\fi
 \else\onlydmatherr@\postdisplayspace\fi}
\def\PrenewCDSpace#1{\global\prenewCDskip@#1\relax}
\def\PostnewCDSpace#1{\global\postnewCDskip@#1\relax}
\def\newCD#1\endnewCD{%
 \outnewCD@\cgaps\outnewCD@\rgaps\outnewCD@\Cgaps\outnewCD@\Rgaps
 \prenewCD@#1\endnewCD
 \advance\abovedisplayskip\prenewCDskip@
 \advance\abovedisplayshortskip\prenewCDskip@
 \advance\belowdisplayskip\postnewCDskip@
 \advance\belowdisplayshortskip\postnewCDskip@
 \vcenter{\vskip\prenewCDskip@ \Let@ \colcount@\@ne \rowcount@\z@
  \everycr{%
   \noalign{%
    \ifnum\rowcount@=\Rowcount@
    \else
     \global\nointerlineskip
     \getrgap@\rowcount@ \vskip\getdim@
     \global\advance\rowcount@\@ne \global\colcount@\@ne
    \fi}}%
  \tabskip\z@
  \halign{&\global\xoff@\z@ \global\yoff@\z@
   \getcgap@\colcount@ \hskip\getdim@
   \hfil\vrule height10\p@ width\z@ depth\z@
   $\m@th\displaystyle{##}$\hfil
   \global\advance\colcount@\@ne\cr
   #1\crcr}\vskip\postnewCDskip@}%
 \prenewCDskip@\z@\postnewCDskip@\z@
 \def\getcgap@##1{\ifcase##1\or\getdim@\z@\else\getdim@\standardcgap\fi}%
 \def\getrgap@##1{\ifcase##1\getdim@\z@\else\getdim@\standardrgap\fi}%
 \let\Width@\relax\let\Height@\relax\let\Depth@\relax\let\Rowheight@\relax
 \let\Rowdepth@\relax\let\Colwdith@\relax
}
\catcode`\@=\active
\hsize 30pc
\vsize 47pc
\def\nmb#1#2{#2}         
\def\cit#1#2{\ifx#1!\cite{#2}\else#2\fi} 
\def\idx{}               
\def\ign#1{}             
\redefine\o{\circ}
\define\X{\frak X}
\define\al{\alpha}
\define\be{\beta}
\define\ga{\gamma}

\define\et{\eta}

\define\la{\lambda}

\define\ch{\chi}
\define\ps{\psi}

\define\Ga{\Gamma}

\define\La{\Lambda}

\define\Ph{\Phi}
\define\Ps{\Psi}
\define\Om{\Omega}
\redefine\i{^{-1}}
\define\row#1#2#3{#1_{#2},\ldots,#1_{#3}}
\define\x{\times}
\define\Id{\operatorname{Id}}
\define\on{\operatorname}
\define\var{\quad}
\define\pr{\operatorname{pr}}
\def\today{\ifcase\month\or
 January\or February\or March\or April\or May\or June\or
 July\or August\or September\or October\or November\or December\fi
 \space\number\day, \number\year}
\topmatter
\title  No slices on the space of generalized connections
\endtitle
\author  Peter W. Michor \\ Hermann Schichl \endauthor
\affil
Institut f\"ur Mathematik, Universit\"at Wien,\\
Strudlhofgasse 4, A-1090 Wien, Austria.\\
Erwin Schr\"odinger Institut f\"ur Mathematische Physik,
Pasteurgasse 6/7, A-1090 Wien, Austria
\endaffil
\address
Institut f\"ur Mathematik, Universit\"at Wien,
Strudlhofgasse 4, A-1090 Wien, Austria
\endaddress
\email Peter.Michor\@esi.ac.at, Hermann.Schichl\@esi.ac.at \endemail
\thanks P.W.M\. was supported by `Fonds zur
F\"orderung der wissenschaftlichen                    
Forschung, Projekt P~10037~PHY'.
\endthanks
\abstract 
On a fiber bundle without structure group the action of the gauge 
group (the group of all fiber respecting diffeomorphisms) on the 
space of (generalized) connections is shown not to admit slices.
\endabstract
\endtopmatter

\document


\subhead \nmb.{1}. Introduction \endsubhead

In modern mathematics and physics actions of Lie groups
on manifolds and the resulting orbit spaces (moduli spaces) are of great
interest. For example, the moduli space of principal connections on a 
principal fiber bundle modulo the group of principal bundle automorphisms is
the proper configuration space for Yang--Mills field theory (as e.g.\ outlined
in \cit!{3}, \cit!{18}, and
\cit!{14}). The structure of these orbit spaces usually
is quite complicated, but sometimes it can be shown that they are stratified
into smooth manifolds. This is usually done by proving a, so called, slice
theorem for the group action.
Also, very recent research in theoretical physics is connected to moduli
spaces: e.g.\ invariance of Euler numbers of moduli spaces of instantons
on 4--manifolds \cit!{19}, moduli spaces of parabolic
Higgs bundles, which are connected to Higgs fields
\cit!{9}, \cit!{13}.
In algebraic topology moduli spaces play an important role, either, 
\cit!{8}, \cit!{16}, \cit!{17}, and also the
definition of the famous Donaldson polynomials involves moduli spaces
(\cit!{2}).

The result presented in this paper is connected to a slice theorem for
the orbit space of connections on a {\it principal\/}
fiber bundle modulo the gauge group,
proved by \cit!{6}.
The situation considered in this paper is a generalization of that.

For a general survey on slice theorems and slices see
\cit!{4}, where a slice theorem for the space of
solutions of Einstein's equations modulo the diffeomorphism group is proved.

The non-existence of the slice theorem in the case of (generalized) connections
on a fiber bundle modulo the gauge group is connected to the fact 
that the right action of $\on{Diff}(S^1)$ on $C^\infty(S^1,\Bbb R)$ 
by composition admits in general no slices, except when restricted to the
space of functions which have finite codimension at all critical points 
\cit!{1}. For more information on (generalized) see \cit!{10}, 
\cit!{11}, or section 9 of \cit!{5}.

\subhead \nmb.{2}. Definitions \endsubhead
A \idx{\it (fiber) bundle}
$(p:E\to M,S)$ consists of smooth finite dimensional manifolds
$E$, $M$, $S$ and a smooth mapping $p: E\to M$. Furthermore, each
$x\in M$ has an open neighborhood $U$ such that $E \mid U :=
p^{-1}(U)$  is diffeomorphic to $U \times S$ via a fiber
respecting diffeomorphism:
$$\Cgaps{.8}\newCD
E \mid U @()\L{\ps}@(2,0) @()\L{p}@(1,-1) & &  
     U \times S @()\l{\pr_1}@(-1,-1)\\
& U. &
\endnewCD$$
In the following we assume that $M$ and $S$, hence $E$ are compact.

We
consider the fiber linear tangent mapping $Tp:TE\to TM$
and its kernel $\ker\,Tp =: VE$, which is called the
\idx{\it vertical bundle\ign{ of a fiber bundle}} of $E$. It is a 
locally splitting vector subbundle of the tangent bundle $TE$.

A \idx{\it connection\ign{ on a fiber bundle}} on the fiber bundle
$(p:E\to M,S)$ is a vector valued 1-form $\Ph \in \Om^1(E;VE)$ with
values in the vertical bundle  $VE$ such that $\Ph\o\Ph=\Ph$ and
$\text{im}\Ph = VE$; so $\Ph$ is just a projection $TE\to VE$.

The kernel $\ker\Ph$ is a sub
vector bundle of $TE$, it is called the space of 
\idx{\it horizontal vectors\ign{ of a fiber bundle}\/} 
or the \idx{\it horizontal bundle}, and it is denoted by 
$HE$. Clearly, $TE=HE\oplus VE$ and $T_uE=H_uE\oplus V_uE$ for $u\in E$.

\comment
Now we consider the mapping $(Tp,\pi_E):TE\to  TM \times_ME$. Then by 
definition 
$(Tp,\pi_E)^{-1}(0_{p(u)},u) = V_uE$, so
$(Tp,\pi_E)\mid HE: HE \to  TM \times_ME$ is a fiber linear
isomorphism, which may be checked in a chart. Its inverse is denoted by 
$$C:= ((Tp,\pi_E)\mid HE)^{-1}: TM \times_M E \to  HE
\hookrightarrow TE.$$   
So $C: TM\times_ME \to  TE$ is fiber linear over $E$ and a right
inverse for $(Tp,\pi_E)$. $C$ is called the 
\idx{\it horizontal lift} associated to the connection $\Ph$. 

Note the formula $\Ph(\xi_u)  = \xi_u - C(Tp.\xi_u,u)$ for
$\xi_u \in T_uE$.   So we can equally well describe a connection
$\Ph$ by specifying $C$. Then we call $\Ph$ 
\idx{\it vertical projection} 
and $\ch := \text{id}_{TE} - \Ph = C\o (Tp,\pi_E)$ will be called  
\idx{\it horizontal projection}.
\endcomment

 If $\Ph: TE \to VE$ is a 
connection on the bundle $(p:E\to M,S)$, then the 
\idx{\it curvature} $R$
of $\Ph$ is given by the Fr\"olicher-Nijenhuis bracket
$$
2R = [\Ph,\Ph] = [\Id - \Ph, \Id - \Ph] \in \Om^2(E;VE).
$$
$R$ is an obstruction against 
involutivity of the horizontal subbundle in the following sense:
If the curvature $R$ vanishes, then horizontal vector 
fields on $E$ also have a horizontal Lie bracket.
\comment
Note that for vector 
fields $\xi,\et\in\X(M)$ and their horizontal lifts 
$C\xi,C\et\in\X(E)$ we have  $R(C\xi,C\et)=[C\xi,C\et]-C([\xi,\et])$. 
\endcomment
Furthermore, we have 
the \idx{\it Bianchi identity} $[\Ph,R] = 0$ by the graded Jacobi 
identity for the Fr\"olicher-Nijenhuis bracket.

\subhead \nmb.{3}. Local description \endsubhead
 Let $\Ph$ be a connection on
$(p:E\to M,S)$. Let us fix a fiber bundle atlas $(U_\al)$ with
transition functions $(\ps_{\al\be})$, and let us consider the
connection $((\ps_\al)^{-1})^*\Ph \in \Om^1(U_\al\times S;
U_\al\times TS)$, which may be written in the form 
$$(((\ps_\al)^{-1})^*\Ph)(\xi_x,\et_y) =: - \Ga^\al(\xi_x,y)  + \et_y
\text{ for }\xi_x \in T_xU_\al \text{ and } \et_y \in T_yS,$$
since it reproduces vertical vectors. The $\Ga^\al$ are given by
$$(0_x,\Ga^\al(\xi_x,y)) := - T(\ps_\al).\Ph.T(\ps_\al)^{-1}.(\xi_x,0_y).$$
We consider $\Ga^\al$ as an element of the space $\Om^1(U_\al;\X(S))$, a
1-form on $U^\al$ with values in the Lie
algebra $\X(S)$ of all vector fields on the standard fiber.
The $\Ga^\al$ are called the \idx{\it Christoffel forms\/} of the
connection $\Ph$ with respect to the bundle atlas $(U_\al,\ps_\al)$.

The transformation law for the Christoffel
forms is 
$$
T_y(\ps_{\al\be}(x,\var)).\Ga^\be(\xi_x,y) =
     \Ga^\al(\xi_x,\ps_{\al\be}(x,y)) - T_x(\ps_{\al\be}(\var,y)).\xi_x.
$$

The curvature $R$ of $\Ph$ satisfies 
$$
(\ps_\al^{-1})^*R = d\Ga^\al + \tfrac12[\Ga^\al,\Ga^\al]_{\X(S)}^\wedge.
$$
Here $d\Ga^\al$ is the exterior derivative of the 1-form 
$\Ga^\al\in\Om^1(U_\al,\X(S))$ with values in the convenient vector 
space $\X(S)$. 

\subheading{\nmb.{4}} The \idx{\it gauge group\ign{ of a fiber 
bundle}\/} $\on{Gau}(E)$ 
of the finite dimensional fiber bundle 
$(p:E\to M,S)$ with compact standard fiber $S$ is, 
by definition, the group of all fiber respecting diffeomorphisms 
$$\CD
E @>f>> E\\
@VpVV @VVpV\\
M@>\Id>> M.
\endCD$$
The gauge group acts on the space of connections by 
$\Ph\mapsto f^*\Ph=Tf\i.\Ph.Tf$. By naturality of the 
Fr\"olicher-Nijenhuis bracket for the curvatures, we have 
$$
R^{f^*\Ph}=\tfrac12[f^*\Ph,f^*\Ph]=\tfrac12f^*[\Ph,\Ph]=f^*R^\Ph = 
Tf\i.R^\Ph.\La^2Tf.
$$ 
Now it is very easy to describe the infinitesimal action. Let $X$ be 
a vertical vector field with compact support on $E$, and consider its 
global flow $\on{Fl}^X_t$. 

Then we have 
$\frac d{dt}|_0(\on{Fl}^X_t)^*\Ph=\Cal L_X\Ph=[X,\Ph]$, the 
Fr\"olicher Nijenhuis bracket, by \cit!{5}. 
The tangent space of $\on{Conn}(E)$ at $\Ph$ is 
the space $T_{\Ph}\on{Conn}(E)=\{\Ps\in\Om^1(E;TE):\Ps|VE=0\}$. The 
"infinitesimal orbit" at $\Ph$ in $T_\Ph\on{Conn}(E)$ is 
$\{[X,\Ph]:X\in C^\infty_c(E\gets VE)\}$.

The isotropy subgroup of a connection $\Ph$ is 
$\{f\in\on{Gau}(E):f^*\Ph=\Ph\}$. Clearly, this is just the 
group of all those $f$ which respect the horizontal bundle 
$HE=\on{ker}\Ph$. It is in general not compact and infinite dimensional.
The most interesting object is of course 
the orbit space $\on{Conn}(E)/\on{Gau}(E)$.

\subhead\nmb.{5}. Slices \endsubhead 
Let $\Cal M$ be a smooth manifold, $G$ a Lie group,  
$G\x \Cal M\to \Cal M$ a smooth action, $x\in \Cal M$, and let  
$G_x=\{g\in G:g.x=x\}$ denote the isotropy group at $x$.  
A contractible subset $S\subseteq \Cal M$ is called a 
\idx{\it slice\/} at $x$, if it contains $x$ and satisfies
\roster
\item"(\nmb:{1})" If $g\in G_x$ then $g.S=S$.
\item"(\nmb:{2})" If $g\in G$ with $g.S\cap S\ne\emptyset$ then
	$g\in G_x$.
\item"(\nmb:{3})" There exists a local section $\chi:G/G_x\to G$ defined
	on a neighborhood $V$ of the identity coset such that the
	mapping $F:V\times S\to\Cal M$, defined by $F(v,s):=\chi(v).s$
	is a homeomorphism onto a neighborhood of $x$.
\endroster
\comment
This is a local version of the definition in finite dimensions, which 
is too narrow for the infinite dimensional situation. However, in finite
dimensions the definition above is equivalent to the usual one:
In the case when $M$ and $G$ are finite dimensional,
a subset $S\subseteq M$ is called a
\idx{\it slice} at $x$, if there is a $G$-invariant open neighborhood  
$U$ of the orbit $G.x$ and a smooth equivariant retraction $r:U \to G.x$ 
such that $S=r\i (x)$. 
In the general case we have the following properties:
\endcomment
We have the following additional properties
\roster
\item"(\nmb:{4})" For 
	$y\in F(V\x S)\cap S$ we get $G_y\subset G_x$, by 
     \therosteritem{\nmb|{2}}.
\item"(\nmb:{5})" For $y\in F(V\x S)$ the isotropy group $G_y$ is 
       conjugate to a subgroup of $G_x$, by \therosteritem{\nmb|{3}} and 
       \therosteritem{\nmb|{4}}.
\endroster
 
\comment
\proclaim{\nmb?{6}. Counter-example}
\cit!{1}{}, \cit!{15}{, 6.2}. 
The right action of $\on{Diff}(S^1)$ on  
$C^\infty(S^1,\Bbb R)$ does not admit slices.
\endproclaim

Let $h(t):S^1=(\Bbb R\mod 1)\to\Bbb R$ be a smooth bump function with 
$h(t) = 0$ for $t\notin[0,\tfrac14]$ and $h(t)>0$ for  
$t\in(0,\tfrac14)$. 
Then put $h_n(t)=\frac1{4^n}h(4^n(t-(1-\frac1{4^n})/3))$ which is  
is nonzero in the interval 
$\bigl((1-\frac1{4^n})/3,(1-\frac1{4^{n+1}})3\bigr)$, and consider  
$$ 
  f_N(t) = \sum_{n=0}^N h_n(t) e^{-\frac1{(t-\frac13)^2}},\qquad 
  f(t) = \sum_{n=0}^\infty h_n(t) e^{-\frac1{(t-\frac13)^2}}. 
$$ 
Then $f\ge 0$ is a smooth function which in $(0,\frac13)$  
has zeros exactly at 
$t=\frac{1-\frac1{4^n}}3$ and which is 0 for $t\notin (0,\frac13)$.  
In every neighborhood of $f$ lies a function $f_N$ 
which has only finitely many of the zeros of $f$  
and is identically zero in the interval 
$[(1-\frac1{4^{N+1}}/3,1/3]$.  
All diffeomorphisms in the isotropy subgroup of $f$ are also contained in 
the isotropy subgroup of $f_N$, but the latter group contains additionally 
all diffeomorphisms of $S^1$ which have support only on 
$[(1-\frac1{4^{N+1}}/3,1/3]$. This contradicts \nmb!{5.5}. 
\endcomment 

\proclaim{\nmb.{6}. Counter-example}
\cit!{15}{, 6.7}. 
The action of the gauge group $\on{Gau}(E)$ on  
$\on{Conn}(E)$ does not admit slices, for $\dim M\ge2$.
\endproclaim

We will construct locally a connection, which satisfies that in any 
neighborhood there exist connections which have a bigger isotropy 
subgroup.  
Let $n=\dim S$, and let $h:\Bbb R^n\to\Bbb R$ be a smooth nonnegative bump function, 
which satisfies $\on{carr} h=\{s\in\Bbb R^n|\;\|s-s_0\|<1\}$. 
Put $h_r(s):=r h(s_0+\tfrac1r(s-s_0))$,  
then $\on{carr} h_r = \{s\in\Bbb R^n|\;\|s-s_0\|<r\}$. 
Then set $h_r^{s_1}(s):= h_r(s-(s_1-s_0))$ 
which implies $\on{carr} h_r^{s_1}=\{s\in\Bbb R^n|\;\|s-s_1\|<r\}$. 
Using these functions, we can define new functions $f_k$ for  
$k\in\Bbb N$ as 
$$ 
    f_k(s) = h_{\|z\|/{2^k}}^{s_k}(s), 
$$ 
where $z:=\tfrac{s_\infty-s_0}3$ for some $s_\infty\in\Bbb R^n$ and 
$s_k := s_0 + z(2\sum_{\l=0}^k \tfrac1{2^{\l}} - 1 - \tfrac1{2^k})$. 
Further define
$$ 
f^N(s) :=e^{-\frac1{\|s-s_\infty\|^2}}\sum_{k=0}^N \tfrac1{4^k}f_k(s),\quad  
f(s) :=\lim_{N\to\infty} f^N(s). 
$$ 
The functions $f^N$ and $f$ are smooth, respectively, since all the 
functions $f_k$ are smooth, and on every point $s$ at most one summand is 
nonzero.  
$\on{carr}f^N=\bigcup_{k=0}^N
     \{s\in\Bbb R^n|\;\|s-s_k\|<\tfrac1{2^k}\|z\|\}$,  
$\on{carr}f=\bigcup_{k=0}^\infty
\{s\in\Bbb R^n|\;\|s-s_k\|<\tfrac1{2^k}\|z\|\}$, $f^N$  
and $f$ vanish in all derivatives in all $x_k$, and $f$ vanishes in all 
derivatives in $s_\infty$.  
 
Let $\ps:E|{U}\to U\times S$ be a fiber bundle chart of  
$E$ with a chart $u:U @>\cong>> \Bbb R^m$ on $M$, and let  
$v:V @>\cong>> \Bbb R^n$ be a chart on $S$. 
Choose $g\in C^\infty_c(M,\Bbb R)$  with 
$\emptyset\ne\on{supp}(g)\subset U$ and   
$dg\wedge du^1\ne 0$ on an open dense subset of $\on{supp}(g)$. 
Then we can define  
a Christoffel form as in \nmb!{3} by 
$$ 
    \Ga := g\,du^1\otimes f(v)\partial_{v^1} \in\Om^1(U,\X(S)). 
$$ 
This defines a connection $\Ph$ on $E|{U}$ which can be extended to a 
connection $\Ph$ on $E$ by the following method.  
Take smooth functions $k_1,k_2\ge0$ on 
$M$ satisfying $k_1+k_2= 1$,  $k_1=1$ on $\on{supp}(g)$, and  
$\on{supp}(k_1)\subset U$, and take an arbitrary connection $\Ph'$ on 
$E$, and set  
$\Ph = k_1\Ph^\Ga+k_2\Ph'$, where $\Ph^\Ga$ denotes the connection 
which is  
induced locally by $\Ga$.  
In any neighborhood of $\Ph$ there exists a connection 
$\Ph^N$ defined by 
$$ 
    \Ga^N := g\,du^1\otimes f^N(s)\partial_{v_1} \in\Om^1(U,\X(S)), 
$$ 
and extended like $\Ph$. 
 
{\bf Claim:} There is no slice at $\Ph$.\newline 
{\it Proof:} We have to consider the isotropy subgroups of $\Ph$ and $\Ph^N$. 
Since the connections $\Ph$ and  
$\Ph^N$ coincide outside of $U$, we may investigate them locally on  
$W=\{u:k_1(u)=1\}\subset U$.  
The curvature of $\Ph$ is given locally on $W$ by \nmb!{3} as
$$ 
R_U:=d\Ga-\tfrac12[\Ga,\Ga]^{\X(S)}_\wedge =  
dg\wedge du^1\otimes f(v)\partial_{v^1} - 0. \tag{\nmb:{1}} 
$$ 
For every element of the gauge group $\on{Gau}(E)$ which is in the  
isotropy group $\on{Gau}(E)_{\Ph}$ the local representative over $W$ 
which looks like 
$\tilde\ga:(u,v)\mapsto (u,\ga(u,v))$ by \nmb!{3} satisfies 
$$\align 
T_v(\ga(u,\var)).\Ga(\xi_u,v) &= 
     \Ga(\xi_u,\ga(u,v)) - T_u(\ga(\var,v)).\xi_u, \tag{\nmb:{2}}\\ 
g(u)du^1\otimes f(v)\sum_i  
     \frac{\partial \ga^1}{\partial v^i}\partial_{v^i} &= 
     g(u)du^1\otimes f(\ga(u,v))\partial_{v^1} - 
     \sum_{i,j}   
     \frac{\partial \ga^i}{\partial u^j}du^j\otimes\partial_{v^i}. 
\endalign$$ 
Comparing the coefficients of $du^j\otimes\partial_{v^i}$ we get the
following equations for $\ga$ over $W$.
$$\align
\frac{\partial\ga^i}{\partial u^j} &= 0 \quad\text{for $(i,j)\ne(1,1)$},\\
g(u)f(v)\frac{\partial\ga^1}{\partial v^1} &=
	g(u)f(\ga(u,v))-\frac{\partial\ga^1}{\partial u^1}. \tag{\nmb:{3}}
\endalign$$
Considering next the transformation $\tilde\ga^*R_U=R_U$ of the 
curvature \nmb!{3} we get 
$$\align 
T_v(\ga(u,\var)).R_U(\xi_u,\eta_u,v) &= R_U(\xi_u,\eta_u,\ga(u,v)), \\
dg\wedge du^1\otimes f(v)\sum_i
     \frac{\partial \ga^1}{\partial v^i}\partial_{v^i} &= 
     dg\wedge du^1\otimes f(\ga(u,v))\partial_{v^1}. \tag{\nmb:{4}}
\endalign$$ 
Another comparison of coefficients yields the equations
$$\align
     f(v)\frac{\partial\ga^1}{\partial v^i} &= 0\quad\text{for $i\ne 1$},\\
     f(v)\frac{\partial\ga^1}{\partial v^1} &= f(\ga(u,v)), \tag{\nmb:{5}}
\endalign$$
whenever $dg\wedge du^1\ne 0$, but this is true on an open dense subset
of $\on{supp}(g)$. Finally, putting \therosteritem{\nmb|{5}} into \therosteritem{\nmb|{3}} shows
$$
    \frac{\partial\ga^i}{\partial u^j} = 0 \quad\text{for all $i$, $j$}.
$$
Collecting the results on $\on{supp}(g)$, we see that $\ga$ has to be constant
in all directions of $u$. Furthermore, wherever $f$ is nonzero, 
$\ga^1$ is a function of $v^1$ only
and $\ga$ has to map zero sets of
$f$ to zero sets of $f$.

Replacing $\Ga$ by $\Ga^N$ we get the same results with
$f$ replaced by $f^N$. Since $f=f^N$ wherever $f^N$ is nonzero or $f$ 
vanishes, $\ga$ in the isotropy group of $\Ph$ obeys all these 
equations not only for $f$ but also for $f^N$ on $\on{supp} f^N\cup f\i(0)$.
On $B:=\on{carr} f\backslash\on{carr} f^N$ the gauge transformation $\ga$ 
is a function of $v^1$ only, hence
it cannot leave the zero set of $f^N$ by construction of $f$ and 
$f^N$.
Therefore, $\ga$ obeys all equations for $f^N$ whenever it obeys all 
equations
for $f$. Thus, every gauge transformation in the isotropy subgroup of
$\Ph$ is in the isotropy subgroup of $\Ph^N$.
 
On the other hand, any $\ga$ having support in $B$ changing only in
$v^1$ direction not keeping the zero sets of $f$ invariant defines a
gauge transformation in the isotropy subgroup of $\Ph^N$ which is not
in the isotropy subgroup of $\Ph$.

Therefore, there exists in every neighborhood of $\Ph$ a connection 
$\Ph^N$ whose isotropy subgroup is bigger than the isotropy subgroup 
of $\Ph$. Thus, by property \nmb!{5.5} no slice exists at $\Ph$.
 
\proclaim{\nmb.{7}. Counter-example}
\cit!{15}{, 6.8}. 
The action of the gauge group $\on{Gau}(E)$ on  
$\on{Conn}(E)$ also admits no slices for $\dim M=1$, i.e.\ for 
$M=S^1$.
\endproclaim

The method of \nmb!{6} is not applicable in this situation, since
there is no function $g$ satisfying $dg\wedge du^1\ne0$ on an open and dense
subset of $\on{supp}(g)$. In this case, however, any connection
$\Ph$ on $E$ is flat. Hence, the horizontal bundle is integrable, the 
horizontal foliation induced by $\Ph$ exists and determines $\Ph$. Any
gauge transformation leaving $\Ph$ invariant also has to map leaves of
the horizontal foliation to other leaves of the horizontal foliation.

We shall construct connections $\Ph^{\la'}$ near $\Ph^\la$ such that 
the isotropy groups in $\on{Gau}(E)$ look radically different near 
the identity, contradicting \nmb!{5.5}. 

Let us assume without loss of generality
that $E$ is connected, and then, by replacing $S^1$ by a finite covering, if 
necessary, that the fiber is connected. 
Then there exists a smooth global section $\chi:S^1\to E$.
By \cit!{12}, p.~95, there exists a tubular neighborhood
$\pi:U\subset E\to\on{im}(\ch)$ such that $\pi=\chi\o p| U$ (i.e.\ a 
tubular neighborhood with vertical fibers). This tubular neighborhood then
contains an open thickened sphere bundle with fiber 
$S^1\x \Bbb R^{n-1}$, and since we are only 
interested in gauge transformations near $\Id_E$, which e.g\. keep a 
smaller thickened sphere bundle inside the larger one, we may replace $E$ 
by an $S^1$-bundle. 
By replacing the Klein bottle by a 2-fold covering, we may finally 
assume that the bundle is $\pr_1:S^1\x S^1\to S^1$.

Consider now connections where the horizontal foliation is a 
1-parameter subgroup with slope $\la$ we see that the isotropy group 
equals $S^1$ if $\la$ is irrational, and equals $S^1$ times the 
diffeomorphism group of a closed interval if $\la$ is rational.

\proclaim{\nmb.{8}. Consequences}
For every compact base manifold $M$ and every compact standard fiber $S$, which
are at least one dimensional, there
are connections on the fiber bundle $E$, where the action of $\on{Gau}(E)$ on
$\on{Conn}(E)$ does not admit slices.
\endproclaim

\enddocument